\documentclass[12pt]{article}
\usepackage{amsmath,amsthm,amsfonts,latexsym,epsfig}
\newtheorem{theorem}{Theorem}[section]

\newtheorem{proposition}[theorem]{Proposition}
\newcommand{\QED}{{\hfill\large $\Box$}}
\title{Dyck paths and a bijection for multisets of hook numbers}
\author{Ian Goulden and Alexander Yong}
\date{\today}
\begin{document}
\maketitle
 
\begin{abstract}
\noindent
We give a bijective proof of a conjecture of Regev and Vershik~
[\ref{RV}] on the equality of two multisets of hook numbers of
certain skew-Young diagrams.
The bijection proves a result that is
stronger and more symmetric than the original conjecture,
by means of a construction
involving Dyck paths, a particular type of lattice path.
\end{abstract}
 
\section{Introduction}
Let $n,k$ be positive integers, and $\alpha =(\alpha_1,\ldots ,\alpha_k )$
be a partition with at most $k$ parts, each part at most $n$,
so $n\geq\alpha_1\geq\ldots\geq\alpha_k \geq 0$. The {\em Young diagram}
of $\alpha$ is given by
$$D=\{ (i,j)\vert 1\leq i\leq k,1\leq j\leq \alpha_{k-i+1}\},$$
a collection of unit cells, arranged in rows and columns. Here
cell $(i,j)$ appears in row $i$ and column $j$, rows numbered
from bottom to top, and columns numbered from left to right. We
regard translates of the diagram in the plane
as equivalent, and generally place
the bottom-left cell at $(1,1)$. (Note, however that this is not the case
for $D$ above when $\alpha_k =0$.) Also let
$$R=\{ (i,j)\vert 1\leq i\leq k,1\leq j\leq n\},$$
$$T=\{ (i,j)\vert 1\leq i\leq k,\alpha_1 -\alpha_i +1
\leq j\leq n+\alpha_1 -\alpha_i\},$$
$$V=\{ (i,j)\vert k+1\leq i\leq 2k,n+\alpha_1-\alpha_{i-k}+1
\leq j\leq n+\alpha_1\},$$
$$SQ=T\cup V,$$
so $R,T,SQ$ are {\em skew diagrams} (in fact, $R$ is also a Young diagram,
the $k\times n$ rectangle).
 
For a skew diagram $G$, let $G^*$ be the skew diagram obtained by rotating
$G$ through 180 degrees. Thus, for example,
$$T^*=\{ (i,j)\vert 1\leq i\leq k,\alpha_{k-i+1} -\alpha_k +1
\leq j\leq n+\alpha_{k-i+1} -\alpha_k\}.$$ 
Also, let $G^{\dagger}$ be
the collection of cells obtained by reflecting $G$ about
a vertical axis.
 
The {\em arm length} $a_G (x)$ of a cell $x$ in a skew diagram $G$
is the number of cells of $G$ in the same row of $x$ and to the right
of $x$; the {\em leg length} $l_G (x)$ of a cell $x$ in a skew diagram $G$
is the number of cells of $G$ in the same column and below.
The {\em coleg} length of a cell $x$ in a skew diagram is the number of
cells
in the same column and above.
The {\em hook length} $h_G (x)$ is given by $h_G (x)=a_G (x)+l_G (x)+1$.
If $E$ is a subset of the cells of $G$, then $AL_G (E)$ is the
multiset $\{ (a_G (x),l_G (x))\vert x\in E\}$, and $H_G (E)$ is
the multiset $\{h_G (x)\vert x\in E\}$.
When there is no ambiguity, we write $H_G (G)$ as $H(G)$,
and $AL_G (G)$ as $AL(G)$.
 
For example, the skew diagrams $D,R,SQ$ are illustrated in Figure 1
for the case $n=6,k=4,\alpha =(6,5,3,1)$. For the three cells
labelled $b,c,d$ in Figure 1, we have $a_D (b)=1,l_D(b)=0,
a_{SQ}(c)=4,l_{SQ}(c)=2$ and $a_R(d)=0,l_R(d)=3$.
 
\begin{figure}[h]
\centering
\epsfig{file=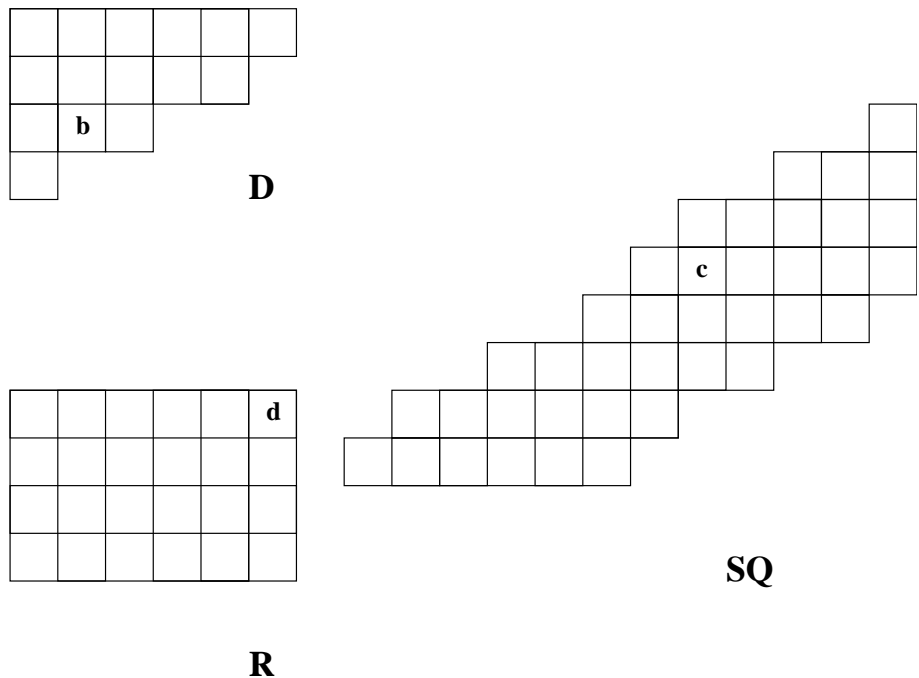,height=9cm}
\caption{$D,R,SQ$ for $n=6,k=4,\alpha =(6,5,3,1)$.}
\end{figure}

Theorem \ref{thmone} below was conjectured by Regev and
Vershik~[\ref{RV}],
and proved by Regev and Zeilberger~[\ref{RZ}],
Janson~[\ref{J}], and Bessenrodt~[\ref{B}] (though only for
the case $n=\alpha_1$ in~[\ref{RZ}]).
 
\begin{theorem}\label{thmone}
For all $n,k,\alpha$,
$$H(SQ)=H(R)\cup H(D)$$
is a multiset identity.
\end{theorem}
\medskip
Regev and Zeilberger note that their proof is not bijective,
and ask for a canonical bijection between the multisets.
Bessenrodt~[\ref{B}] presents such a bijection, deducing it from
a general result about ``removable'' hooks in Young diagrams. In
this paper, we present a different bijection, deducing it from
another general result, the main result of the paper.
It is convenient to keep
arm and leg lengths separately, and thus we prove the following
result, which is obviously a generalization of Theorem \ref{thmone}.
 
\begin{theorem}\label{thmtwo}
For all $n,k,\alpha$,
$$AL(SQ)=AL(R)\cup AL(D)$$
is a multiset identity.
\end{theorem}
 
\medskip
The next result, our main result, is more symmetric and
natural looking than Theorem~\ref{thmtwo}, but it implies
Theorem~\ref{thmtwo}.
Independently, this result has also been obtained by Regev~[\ref{R1}],
and a bijective proof that is different from ours has been given
by Krattenthaler~[\ref{K}].
 
\begin{theorem}\label{thmthr}
For all $n,k,\alpha$,
$$AL(T)=AL(T^*)$$
is a multiset identity.
\end{theorem}
 
\medskip
We delay the proof of Theorem~\ref{thmthr} until the next section, and
proceed now by giving a bijective proof that it implies
Theorem~\ref{thmtwo}. The proof involves partitioning the cells 
of $R$ and $T^*$ into two regions each, and identifying cells in 
various regions of skew diagrams whose pairs of arm and 
leg lengths are immediately
equal. 

\medskip
\noindent{\bf Proof that Theorem~\ref{thmthr} implies
Theorem~\ref{thmtwo}:}
Partition the cells of $R$ into two subsets $R_1$ and $R_2$, given by
$$R_1=\{ (i,j)\vert 1\leq i\leq k,n-\alpha_{k-i+1}+1\leq j\leq n\},$$
$$R_2=\{ (i,j)\vert 1\leq i\leq k,1\leq j\leq n-\alpha_{k-i+1}\},$$
and the cells of $T^*$ into two subsets $T^*_1$ and $T^*_2$, given by
$$T^*_1=\{ (i,j)\vert 1\leq i\leq k,\alpha_{k-i+1}-\alpha_k +1
\leq j\leq n-\alpha_k\},$$
$$T^*_2=\{ (i,j)\vert 1\leq i\leq k,n-\alpha_k +1
\leq j\leq n+\alpha_{k-i+1}-\alpha_k\}.$$
The significance of these regions in this proof is
that $R_1^{\dagger}=T_2^*=V^*=D$ and $R_2^{\dagger}=T_1^*$.
These equalities (using appropriate translations) are immediate from the
definitions of the regions.
See Figure 2 for an illustration of these regions in the
case $n=6,k=4,\alpha =(6,5,3,1)$, and to check visually the above
equalities
in this case.
 
\begin{figure}[h]
\centering
\epsfig{file=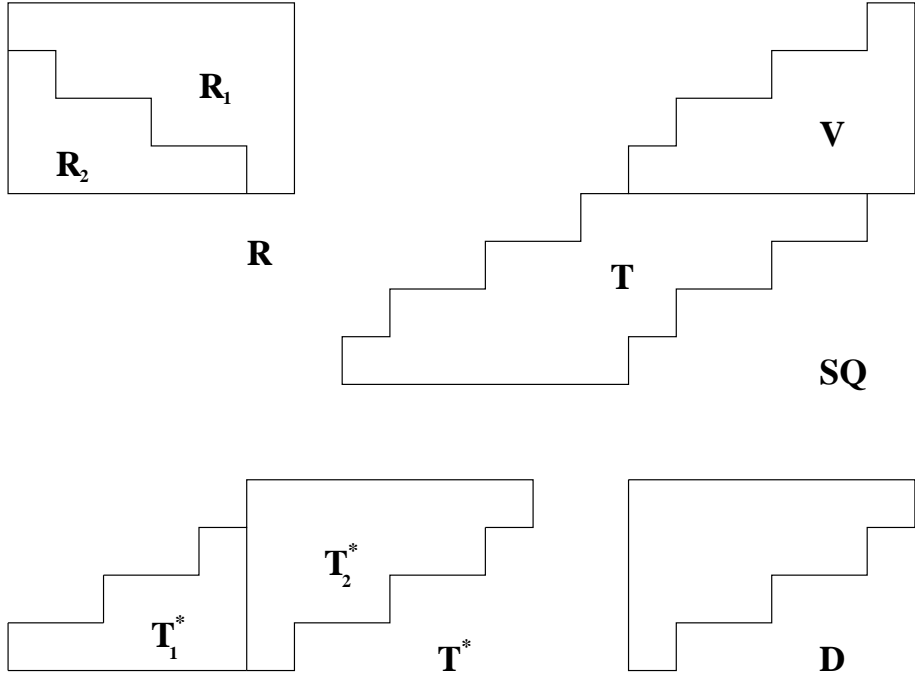,height=9cm}
\caption{Skew shapes for $n=6,k=4,\alpha =(6,5,3,1)$.}
\end{figure}

\medskip
\noindent{\it Bijective identification of $AL_{SQ}(V)$ and $AL_R(R_1 )$:}
Now $V^*=R_1^{\dagger}$, so
the $j$th columns of $V$ and $R$, respectively, have
the same lengths, for each $j=1,\ldots ,\alpha_1$. Furthermore, 
$V$ appears in $SQ$ with cells added below $V$ to extend all columns of
$V$ 
to length $k$. Similarly,
$R_1$ appears in $R$ with cells added below $R_1$ to
extend all columns of $R_1$ 
to length $k$. Thus the arm and leg lengths are equal, for the cells
that are $i$ rows from the topmost entry, in the 
$j$th column from the left most column, of $V$ in $SQ$ and $R_1$ in $R$, 
respectively.
Thus we establish immediately that
\begin{equation}\label{alsqr}
AL_{SQ}(V)=AL_R(R_1 ).
\end{equation}
 
\medskip
\noindent{\it Bijective identification of $AL_{T^*}(T^*_1)$ and
$AL_{R}(R_2)$:}
Now $T_1^*=R_2^{\dagger}$, so the $i$th rows
of $T_1^*$ and $R_2$, respectively, have
the same lengths, for each $i=1,\ldots ,k$ (some of these lengths are zero
when $\alpha_1 =n$). Furthermore, $T_1$ appears in $T$ with cells added
to the right of $T_1$ to extend all rows of $T_1$ to length $n$.
Similarly, $R_2$ appears in $R$ with cells added
to the right of $R_2$ to extend all rows of $R_2$ to length $n$.
Thus, the arm lengths and leg lengths are equal, for
the cells that are $j$ columns from the left most entry, in
the $i$th row from the bottom
row, of $T^*_1$ in $T^*$ and $R_2$ in $R$, respectively.
Thus we establish immediately that
\begin{equation}\label{alrts}
AL_{T^*}(T^*_1)=AL_R(R_2 ).
\end{equation}
 
\medskip
\noindent{\it Bijective identification of $AL_{T^*}(T^*_2 )$ and $AL(D)$:}
Now $T_2^*=D$, and $T^*_2$ appears in $T^*$ with no cells added
to the right nor below, so we establish immediately that
\begin{equation}\label{aldts}
AL_{T^*}(T_2^*)=AL(D).
\end{equation}
\medskip
\noindent{\it The result:} Suppose Theorem~\ref{thmthr} is true. Then,
applying~(\ref{alsqr}), we obtain
\begin{equation}\label{altts}
AL_{SQ}(V)\cup AL(T)=AL_{R}(R_1 )\cup AL(T^*).
\end{equation}
But $AL(T)=AL_{SQ}(T)$, since $T$ appears in $SQ$ with
no cells added to the right nor below.
Also, $AL(T^*)=AL_{T^*}(T_1^*) \cup AL_{T^*}(T_2^*)$, since $T_1^*$ and 
$T_2^*$ partition the cells of $T^*$. Making these substitutions
into~(\ref{altts}) gives
\begin{eqnarray*}
AL_{SQ}(V) \cup AL_{SQ}(T)&=&AL_{R}(R_1 )\cup
AL_{T^*}(T_1^*) \cup AL_{T^*}(T_2^*)\\
&=&AL_{R}(R_1 )\cup AL_{R}(R_2 )\cup AL(D),
\end{eqnarray*}
with the second equality from~(\ref{alrts}) and~(\ref{aldts}).
Now $V$ and $T$ partition
the cells of $SQ$, and $R_1$ and $R_2$ partition the cells of $R$, so
the above result becomes $AL(SQ)=AL(R) \cup AL(D)$, and
we have established Theorem~\ref{thmtwo}.
\QED
 
\medskip
How is this proof bijective? To prove Theorem~\ref{thmthr} bijectively, in
the next section we determine an explicit
bijection $\phi: T \rightarrow T^*$, that preserves arm
and leg lengths (this means
that for each cell $x\in T$ we have $a_T (x)=a_{T^*}(\phi (x) )$ and 
$l_T (x)=l_{T^*}(\phi (x) )$ ). Similarly, to give a bijective proof of
Theorem~\ref{thmtwo}, we must determine an explicit
bijection $\psi : SQ \rightarrow R \cup D$, that preserves arm
and leg lengths.
 
In terms of $\phi$, we now describe such a
bijection $\psi$ that is implicit in the above proof.
First, note that, to establish~(\ref{alsqr}), (\ref{alrts})
and (\ref{aldts}) above, we
have described three simple bijections, and let us call
them $\zeta _1: V\rightarrow R_1$, $\zeta _2: T_1^*\rightarrow R_2$, and
 $\zeta_3: T_2^*\rightarrow D$.
 
\medskip
\noindent {\bf A bijection $\psi$ that establishes Theorem~\ref{thmtwo}:}
For $x\in SQ$, we obtain
 $\psi (x)\in R\cup D$ as follows:
 
\begin{description}
\item For $x\in V$, let $\psi (x)=\zeta_1 (x)$.
\item For $x\in T$,
\begin{itemize}
\item if $\phi (x)\in T_1^*$, let $\psi (x)=\zeta_2 (\phi(x))$,
\item if $\phi (x)\in T_2^*$, let $\psi (x)=\zeta_3 (\phi(x))$.
\end{itemize}
\end{description}
 
\noindent This clearly specifies a bijection $\psi$ of the required type,
giving a bijective proof of Theorem~\ref{thmtwo}.
 
\section{Dyck paths and the bijection}
 
In this section, we determine a bijection $\phi: T \rightarrow T^*$,
that preserves arm and leg lengths, as referred to above at the end
of Section 1. This provides a bijective proof of Theorem~\ref{thmthr}.
 
The bijection is described in terms of a particular type of
lattice path that will be associated with $T$ and $T^*$, called
a Dyck path.
A {\em Dyck path} of length $2k,k\geq 0$, is a sequence $(i,y_i),
i=0,\ldots ,2k$, of lattice points in the plane, in which
$y_0 =y_{2k} =0,y_i\geq 0$, for $i=1,\ldots ,2k-1$, and $y_i -y_{i-1}=+1$
or $-1$, for $i=1,\ldots ,2k$. Equivalently, a Dyck path is
completely specified by its sequence of {\em steps};
if $y_i -y_{i-1} =+1$ then the $i$th step is an up step, and
if $y_i -y_{i-1} =-1$ then the $i$th step is a down step.
The {\em height} of the $i$th step is $y_{i-1}$, for $i=1,\ldots ,2k$.
Since $y_{2k}=0$, then the $2k$ steps consist of $k$ up steps
and $k$ down steps.
We can visualize a Dyck path as a connected path in the plane by
drawing a line segment between the consecutive lattice points in
the path.
 
Let the skew diagrams $T_{[i]}$ and $T_{(i)}$, for $i=1,\ldots ,n$,
be given by
$$T_{[i]} =\{x\in T\vert a_T (x) =i-1\},$$
$$T_{(i)} =\{x\in T\vert a_T (x)\leq i-1\},$$
and define $(T^*)_{[i]}$ and $(T^*)_{(i)}$ in the same way.
%
%
Consider the skew diagram $T_{(i)}$, for each fixed $i=1,\ldots ,n$.
Label the $k$ cells of $T_{[i]}$ in $T_{(i)}$, successively, $x_1,\ldots
,x_k$,
from bottom to top (there is exactly one cell of $T_{[i]}$ in each of
the $k$ rows of $T_{(i)}$). Label the cells of $T_{[0]}$ in $T_{(i)}$,
successively, $z_1,\ldots ,z_k$, from top to bottom (similarly, there is
exactly one cell of $T_{[0]}$ in each of the $k$ rows of $T_{(i)}$). In
the case $i=1$, then each cell of $T_{[0]}$ will have two labels, one
an $x_j$ and the other $z_{k+1-j}$, for some $j=1,\ldots ,k$.
 
Now form a permutation $\sigma_i$ of $x_1,\ldots ,x_k,z_1,\ldots ,z_k$ as
follows: Place the $x$'s and $z$'s from left to right in $\sigma_i$ in
the order that they appear from left to right as labels in the cells
of $T_{(i)}$. For labels in the same column of $T_{(i)}$, order them
with the $x$'s first, followed by the $z$'s; the $x$'s are ordered as
they appear from bottom to top in the same column, and the $z$'s from
bottom to top also.
For example, in the case $n=11,k=9, \alpha =(11,11,9,8,8,6,3,1,0)$, we
illustrate $T_{(3)}$ in Figure 3, with the cells labelled as described
above.
In this case, the permutation $\sigma _3$  is given by
$$\sigma _3 = x_1 x_2 x_3 z_9 z_8 x_4 x_5 z_7 x_6 z_6 z_5 z_4 x_7 x_8 z_3
x_9 z_2 z_1.$$

\begin{figure}[h]
\centering
\epsfig{file=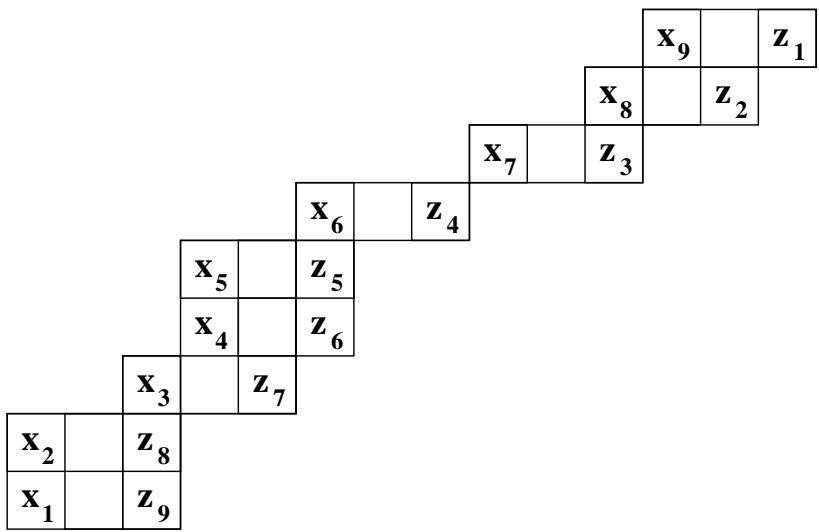,height=7cm}
\caption{$T_{(3)}$ for $n=11,k=9, \alpha =(11,11,9,8,8,6,3,1,0)$.}
\end{figure}

Now let $\rho_i$ be the lattice path starting at $(0,0)$, whose steps
are specified by $\sigma_i$ as follows: the $x_j$'s specify
the up steps (labelled $x_j$),
and the $z_j$'s specify the down steps (labelled $z_j$).
For example the lattice path $\rho_3$ determined from $\sigma_3$ in the
example above is illustrated in Figure 4.
 
\begin{figure}[h]
\centering
\epsfig{file=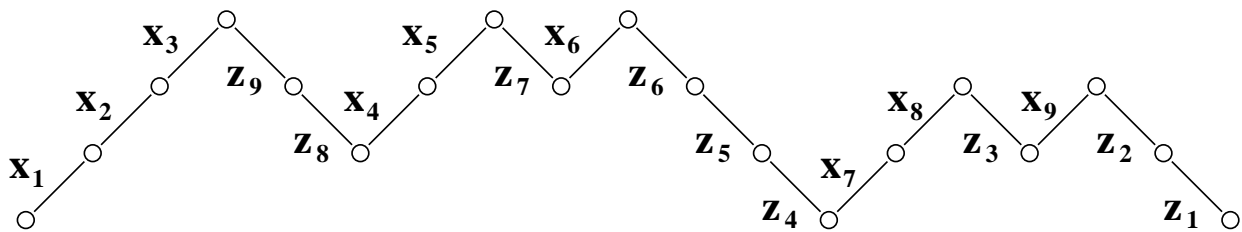,height=3cm}
\caption{The Dyck path $\rho_3$ determined from $\sigma_3$.}
\end{figure}

It is a straightforward induction to prove that the height of the up step
labelled
$x_j$ in $\rho_i$ is equal to the leg length of the cell labelled $x_j$ in
$T_{(i)}$, and that the height of the down step labelled $z_j$ in $\rho_i$
is equal to one more than the coleg length of the cell
labelled $z_j$ in $T_{(i)}$. But since leg and coleg lengths are always
nonnegative, the height of every up step in $\rho_i$ is nonnegative, 
and the height of every down step in $\rho_i$ is positive,
so $\rho_i$ is a Dyck path. For example, the lattice path $\rho_3$
illustrated
in Figure 4 is clearly a Dyck path.
 
Now there is a natural bijection between the  up steps and down steps
in a Dyck path: pair each up step at height $j$ with the first down step
at
height $j+1$ occurring after that up step (there must be such a down step
since the path ends at a vertex with ordinate equal to $0$, and down steps
decrease the value of the ordinate by exactly 1 for each step).
Suppose that the up step labelled $x_j$ is paired with the down step
labelled $z_{P_i(j)}$ in this way, for $j=1,\ldots ,k$. Then $P_i$ is
a bijection on $\{ 1,\ldots ,k\}$, for each fixed $i$. For example, for
the Dyck path illustrated in Figure 4, we have $P_3(1)=4$, $P_3(2)=8$,
$P_3(3)=9$, $P_3(4)=5$, $P_3(5)=7$, $P_3(6)=6$, $P_3(7)=1$, $P_3(8)=3$,
and $P_3(9)=2$.
 
Now rotate $T_{(i)}$, with its cells labelled as above, through 180
degrees,
to obtain $\delta$. Now $\delta =(T_{(i)})^* = (T^*)_{(i)}$, and the cells
of $T_{[0]}$ in $T_{(i)}$, labelled  with $z_j$'s, become
the cells of $(T^*)_{[i]}$ in $\delta$. Moreover, the coleg length 
of a cell labelled
$z_j$ in $T_{(i)}$ equals the leg length of the corresponding cell in
$\delta$,
so
$$l_{T_{(i)}}(x_j )=l_{(T^*)_{(i)}}(z_{P_i (j)}),$$
where, for example, $l_{T_{(i)}}(x_j )$ means the leg length of the cell
labelled $x_j$ in $T_{(i)}$.
Also,
$$a_{T_{(i)}}(x_j )=i-1=a_{(T^*)_{(i)}}(z_{P_i (j)}),$$
since all cells in $T_{[i]}$ and $(T^*)_{[i]}$ have arm length
equal to $i-1$, for each fixed $i$.
But $T_{(i)}$ appears in $T$ with no cells added to the right nor
below, so $l_{T_{(i)}}(x_j )=l_T(x_j )$ and $a_{T_{(i)}}(x_j )=a_T(x_j )$. 
Similarly, $l_{(T^*)_{(i)}}(z_{P_i (j)})=
l_{T^*}(z_{P_i (j)})$ and $a_{(T^*)_{(i)}}(z_{P_i (j)})= a_{T^*}(z_{P_i
(j)})$.
Thus, putting these equalities together, we have
\begin{equation}\label{lastpiece}
l_T(x_j )=l_{T^*}(z_{P_i (j)}),\;\;\;a_T(x_j )=a_{T^*}(z_{P_i (j)}).
\end{equation}

\medskip
\noindent{\bf Proof of Theorem 1.3:} This follows from Lemma 2.1
immediately,
These equations imply that the mapping from the cell labelled $x_j$ in
 $T$ to the cell labelled $z_{P_i (j)}$ in $T^*$, for each $i=1,\ldots
,n$,
is arm and leg length preserving, so we have found the bijection
$\phi$ that we require, as stated below. \QED

\medskip
\noindent {\bf A bijection $\phi$ that establishes Theorem~\ref{thmthr}:}
For $w\in T$, we obtain $\phi (w)\in T^*$ as follows. Each $w$ is
contained in $T_{[i]}$ for some unique $i=1,\ldots ,n$. If $w$ has
label $x_j$ in $T_{(i)}$, then $\phi (w)$ is the cell with
label $z_{P_i (j)}$ in $(T^*)_{(i)}$.
 
\medskip
\noindent This clearly specifies a bijection, that is arm and leg length 
preserving from~(\ref{lastpiece}), giving a bijective proof
of Theorem~\ref{thmthr}.

\section{The projective case}

A refinement of Theorem~\ref{thmtwo} has been given by Regev~[\ref{R2}],
in which the partition $\alpha$ has a special form.
In order to state this result, we require some adaptations of the notation
in Section 1. Let $n=k+1$, and $\alpha$ have the form $\alpha=
(\lambda_1,\ldots , \lambda_m \vert \lambda_1 -1,\ldots ,\lambda_m -1)$,
in {\em Frobenius notation}, where $k\geq \lambda_1 > \ldots > \lambda_m
> 0$, so $\lambda =(\lambda_1, \ldots , \lambda_m )$ is a partition
with $m$ distinct parts. This means that $D$, the Young diagram of $\alpha$,
has exactly $m$ cells on the  (top-left to bottom-right) {\em diagonal}, given
by the cells $(k+1-j,j)$, for $j=1,\ldots ,m$, with $\lambda_j$ cells to the
right of the $j$th of these cells in row $k+1-j$, and $\lambda_j-1$ cells
below this cell in column $j$. Let $\mathcal{B}$ consist of all partitions
 $\alpha$ of this form, for any $m\geq 0, k\geq 1$ (e.g., $R$ is
the Young diagram of a partition in $\mathcal{B}$, with
 $m=k$ and $\lambda_j=k+1-j$, for $j=1,\ldots ,k$).

For a Young diagram $G$, let $p(G)$ consist of the cells of $G$ on or below
the diagonal (as described above), and let $q(G)$ consist of the cells
of $G$ strictly above the diagonal. For a skew diagram, extend this
notation by describing the diagonal: for $T, T_{[i]}, T_{(i)},
 SQ$, where $n=k+1$, and
 $\alpha\in\mathcal{B}$, the
diagonal consists of the cells $(k+1-j,\alpha_1 +j)$, for $j=1, \ldots, k-m$; for
 $T^*$, the diagonal consists of the cells $(k+1-j,k+1-\alpha_k+j)$, for $j=1,
\ldots ,m$.
For example, the skew diagrams $D, R, SQ, T$ are illustrated in Figure XXX
for the case $k=5, m=2, \alpha =(5,4,2,1)$, corresponding
to $\lambda =(4,2)$. In each of these skew diagrams, there is a thick line
extending from top left to bottom right, which partitions the diagram $G$ into
the cells of $p(G)$, below and to the left of the line, and the cells
of $q(G)$, above and to the right of the line.

The following result has been given by Regev~[\ref{R2}], whose proof
is not bijective. A bijective proof has been given by Krattenthaler~[\ref{K}].

\begin{theorem}\label{regnew}
For all $k,m$ and $\alpha\in\mathcal{B}$,
$$AL(p(SQ))=AL(p(R))\cup AL(q(D))$$
is a multiset identity.
\end{theorem}

In order to prove Theorem~\ref{regnew}, we first
note that

\begin{equation}\label{TSQ}
AL(p(SQ))=AL(p(T)),
\end{equation}

so we shall work with $T$ on the left hand side of the result, instead of $SQ$.
For each $i=1,\ldots ,k+1$, let $u$ be the smallest row index among the elements
of $T_{[i]}$ above the diagonal of $T$. Let $T^i$ be the skew diagram obtained
from $T$ by shifting rows $u, u+1 , \ldots , k$ to the right, where necessary,
so that the right most of the $k+1$ cells in each of these rows occurs in
column $\alpha_1 +k+1$. (If no element of $T_{[i]}$ is above the diagonal
of $T$, then we define $T^i=T$.) The diagonals of $T^i$ and ${T^i}^*$ are
the same as for $T$ and $T^*$, respectively.
For example, the skew diagrams $D, R, T,
T^i, {T^i}^*$ are illustrated in Figure XXX for the case $k=12, m=6,
\alpha = (12,11,11,9,8,8,6,4,3,3,1)$,
with $i=5$. In each of these skew diagrams, there is again a
thick line partitioning the cells into those given by $p$ and $q$, and there
is a dot in every cell with arm length equal to $i-1 =4$.

We require the following technical result about the row index $u$, chosen
above for each $i$.

\begin{proposition}\label{techprop}
Let $\alpha\in\mathcal{B}$, with the diagonal of length $m$, and with
 $\alpha_1\leq k+1$. Let $u$ be the smallest row index among the elements
of $T_{[i]}$ above the diagonal of $T$. Then

\begin{enumerate}
\item
$u-\alpha_u >i-1\;\;$ and $\;\; u-1-\alpha_{u-1}\leq i-1$,
\item
$u > m$,
\item
$\alpha_{u-i}\geq u\;\;$ and $\;\;\alpha_{u-i+1}\leq u$,
\item
$\alpha_u +i\leq \alpha_{u-i}\;\;$ and $\;\;\alpha_{u-1}+i\geq \alpha_{u-i+1}$,
\end{enumerate}

\end{proposition}

\noindent{\bf Proof:}
In the row of $T$ with index $a$, for $a=1,\ldots ,k$, the diagonal cell
is in column $\alpha_1 +k+1-a$, the right most element is in 
column $\alpha_1 +k+1-\alpha_a$, and the unique element of $T_{[i]}$ is
therefore in column $\alpha_1 +k+1-\alpha_a -(i-1)$. This means that
the element of $T_{[i]}$ in row $a$ is above the diagonal of $T$ exactly
when $\alpha_1 +k+1-\alpha_a -(i-1)>\alpha_1 +k+1-a$, or
 $a-\alpha_a>i-1$. Part 1 of the result follows immediately.

>From Part 1, we have $u-\alpha_u >i-1\geq 0$, so $\alpha_u < u$. But, since
 $\alpha\in\mathcal{B}$, then $\alpha_j \geq j$ for $j=1,\ldots ,m$, where
 $m$ is the length of the diagonal of $\alpha$, giving Part 2 of the result.

Now let $u-\alpha_u=c$ and $u-1-\alpha_{u-1}=d$, where $c>i-1\geq d$, from
Part 1. Thus in the Young diagram $D$ of$\alpha$, the right most cell in
row $k+1-u$ is in column $u-c$, and the right most cell in row $k+1-(u-1)$
 is in column $u-1-d$. But $\alpha\in\mathcal{B}$, so symmetry of
 $\mathcal{B}$ implies that the bottom cell in column $u+1$ is in
row $k+1-(u-c)$, and the bottom cell in column $u$ is in
row $k+1-(u-1-d)$. Thus we have $\alpha_{u-c}\geq u+1,\;\alpha_{u+1-c}=
\ldots =\alpha_{u-1-d}=u,\;\alpha_{u-d}<u$, and Result 3 follows
from $c>i-1\geq d$.

Part 4 follows immediately from Parts 1 and 3.
\QED

\medskip
Now we are able to prove Theorem~\ref{TSQ}, using the bijective proof of
Theorem~\ref{thmthr}.

\medskip
\noindent{\bf Proof of Theorem~\ref{TSQ}:}
Let $M_1,M_2,M_3,M_4$ be the multisets of leg lengths of the cells
with arm lengths equal to $i-1$, in $T^i,(T^i)^*,p(T),q(D)$, respectively.
Now, Theorem~\ref{thmthr} applied to skew diagram $T^i$ gives a bijection
between $AL(T^i)$ and $AL((T^i)^*)$, which contains a bijection between
 $M_1$ and $M_2$. 

Now, the elements of $M_1$ can be partitioned into two subsets: $M_{11}$,
corresponding to the cells on or below the diagonal of $T^i$; and $M_{12}$,
corresponding to the cells above the diagonal. Thus the elements
of $M_{11}$ correspond to cells in rows $1,\ldots ,u-1$ of $T^i$, and the
elements of $M_{12}$ correspond to the cells in rows $u,\ldots ,k$. But
 $T$ and $T^i$ differ only in rows $u,\ldots ,k$, so $M_{11}=M_3$. Also,
the right most cell of $T^i$ is in column $k+1+\alpha_1 -\alpha_j$,
for $j=1,\ldots , u-1$. Now let $s$ be chosen so that

\begin{equation}\label{sdef}
\alpha_s\leq i-1\;\;and \;\;\alpha_{s-1}>i-1.
\end{equation}

Then the bottom element of column $k+1+\alpha_1-(i-1)$ in $T^i$ is in
row $s$, so $M_{12}=\{ u-s,\ldots ,k-s\}$, giving

\begin{equation}\label{firstset}
M_1=M_3\cup\{ u-s,\ldots ,k-s\}.
\end{equation}

Similarly, the elements of $M_2$ can be partitioned into three
subsets: $M_{21}$, corresponding to the cells in columns $1,\ldots ,k+1$ of
 ${T^i}^*$; $M_{22}$, corresponding to the cells to
the right of column $k+1$ but on or below the diagonal of ${T^i}^*$;
and $M_{23}$, corresponding to the cells above the diagonal of ${T^i}^*$.
Now, the right most cell in rows  $1,\ldots ,k+1-u$  of ${T^i}^*$ is
in column $k+1$, and the right most cell in row $j$ of ${T^i}^*$ is
in column $k+1+\alpha_{k+1-j}$, for $j=k+2-u,\ldots ,k$. Therefore,
from~(\ref{sdef}), the cells in $M_{21}$ occur in
rows $1,\ldots ,k+1-s$, and the bottom element in each corresponding
column is in row $1$, so $M_{21}=\{ 0,\ldots ,k-s\}$.

Now, let $r$ be the largest row index of the elements of $M_{22}$. Then,
since the diagonal elemens of row $j$ is in column $k+1+k+1-j$, for
 $k=k+2-u,\ldots ,k$, we have

\begin{equation}\label{rdef}
k+1-r+i-1\geq\alpha_{k+1-r}\;\;and\;\;k+1-(r+1)+i-1<\alpha_{k-r},
\end{equation}
and, from Proposition~\ref{techprop}(3), we immediately have $k-r=u-i$, or
 $r=k-u+i$. Also, the bottom element of the columns corresponding to the
cells of $M_{22}$ all occur in row $k+2-u$, from
the second part of Proposition~\ref{techprop}(4). Thus,
$M_{22}=\{(k+2-s)-(k+2-u),\ldots ,(k-u+i)-(k+2-u)\} =\{u-s,\ldots ,i-2\}$.

Finally, the leg lengths of the cells of $M_{23}$ are all the same
in ${T^i}^*$ as in $T^*$, from the first part of Proposition~\ref{techprop}(4).
Thus $M_{23}=M_4$, and we have
$$M_{2}=M_{21}\cup M_{22}\cup M_{23}=M_4\cup\{0,\ldots ,k-s\}\cup\{ u-s,ldots ,
i-2\}.$$
The bijection between $M_1$ and $M_2$ then gives, from~(\ref{firstset}),
$$M_3\cup\{ u-s,\ldots ,k-s\}=M_4\cup\{0,\ldots ,k-s\}\cup\{ u-s,ldots ,i-2\},$$
and we have
$$M_3=M_4\cup\{0,\ldots ,i-2\}.$$
Now, Theorem~\ref{regnew} follows from~(\ref{TSQ}), and the fact that the cells in $p(R)$
 with arm length equal to $i-1$ in $R$ have leg lengths $0,\ldots ,i-2$.
\QED

\section*{Acknowledgements}
 
This work was supported by the Natural Sciences and Engineering
Research Council of Canada, through a grant to IG, and a PGSA
to AY.
 
\section*{References}
 
\begin{enumerate}
\item\label{B}
C. BESSENRODT, On Hooks of Young Diagrams, Annals of Combinatorics
{\bf 2}(1998), 103--110.
\item\label{J}
S. JANSON, Hook lengths in a Skew Young Diagram, Electronic J. Combinatorics
{\bf 4}(1997), R24.
\item\label{K}
C. KRATTENTHALER, Bijections for hook pair identities, Electronic J.
Combinatorics {\bf 7}(2000), R27.
\item\label{R1}
A. REGEV, Generalized hook and content numbers identities, European J.
Combinatorics, (to appear).
\item\label{R2}
A. REGEV, Generalized hook and content numbers identities -- The projective
case, European J. Combinatorics, (to appear).
\item\label{RZ}
A. REGEV and D. ZEILBERGER, Proof of a conjecture on Multisets of Hook
Numbers, Annal of Combinatorics {\bf 1}(1997), 391--394.
\item\label{RV}
A. REGEV and A. VERSHIK, Asymptotics of Young's diagrams and hook numbers,
Electronic J. Combin. {\bf 4}(1997), R22.
\end{enumerate}

\end{document}